\documentclass[12pt]{article}
\usepackage{rotating,epic,eepic}
\usepackage{natbib}
\usepackage{amsmath}
\usepackage{amssymb}
\usepackage{delarray}
\usepackage{amsfonts}
\usepackage{amsthm}
\usepackage[english]{babel}
\usepackage{enumerate}
\usepackage{fancyhdr}

\def\hf{\hat{f}_n}

\newcommand \cA{{\cal A}}

\newcommand \cC{{\cal C}}

\newcommand \cE{{\cal E}}
\newcommand \cF{{\cal F}}

\newcommand \cI{{\cal I}}

\newcommand \cN{{\cal N}}
\newcommand \cO{{\cal O}}
\newcommand \cP{{\cal P}}

\newcommand \cS{{\cal S}}

\newcommand \cX{{\cal X}}

\newcommand{\1}{{\rm 1}\kern-0.24em{\rm I}}
\newcommand{\hfn}{{\hat{f}_n}}

\theoremstyle{plain}
\newtheorem{theo}{Theorem}
\newtheorem{coro}{Corollary}

\newtheorem{pro}{Proposition}
\newtheorem{rem}{Remark}
\newtheorem{defi}{Definition}

\date{}

\begin{document}

\title{{\bf Classification with Minimax Fast Rates for Classes of Bayes Rules with Sparse Representation}}
\author{
  Guillaume Lecu\'e\\
{\it Universit\'e Paris VI } \\
}

\maketitle
\date
\bigskip
\begin{abstract}\noindent
We construct a classifier which attains the rate of convergence
$\log n/n$ under sparsity and margin assumptions. An approach
close to the one met in approximation theory for the estimation of
function is used to obtain this result. The idea is to develop the
Bayes rule in a fundamental system of $L^2([0,1]^d)$ made of
indicator of dyadic sets and to assume that coefficients, equal to
$-1,0 \mbox{ or } 1$, belong to a kind of $L^1-$ball. This
assumption can be seen as a sparsity assumption, in the sense that
the proportion of coefficients non equal to zero decreases as
"frequency" grows. Finally, rates of convergence are obtained by
using an usual trade-off between a bias term and a variance term.
\end{abstract}

\section{Introduction}
Consider a measurable space $(\cX,\cA)$ and $\pi$ a probability
measure on this space. Denote by $D_n=(X_i,Y_i)_{1\leq i \leq n}$
$n$ observations of $(X,Y)$ a random variable with values in
$\cX\times\{-1,1\}$ distributed according to $\pi$. We want to
construct measurable functions which associate a label
$y\in\{-1,1\}$ to each point $x$ of $\cX$, such functions are called
{\emph{prediction rules}}. The quality of a prediction rule $f$ is
given by the value
$$R(f)=\mathbb{P}(f(X)\neq Y)$$ called \emph{misclassification error of $f$}.
It is well known (e.g. \cite{dgl:96}) that there exists an optimal
prediction rule which attains the minimum of $R$ over all measurable
functions with values in $\{-1,1\}$. It is called \emph{Bayes rule}
and defined by
$$f^{*}(x)=\rm{sign}(2\eta(x)-1),$$ where $\eta$ is the \emph{conditional probability function
of $Y=1$ knowing $X$} defined by $$\eta(x)=\mathbb{P}(Y=1\vert
X=x).$$ The value
$$R^{*}=R(f^{*})=\min_{f}R(f)$$ is known as the \emph{Bayes risk}. The aim
of classification is to construct a prediction rule, using the
observations $D_n$, which has a risk as close to $R^*$ as possible.
Such a construction is called a {\emph{classifier}}. Performance of
a classifier $\hat{f}_n$ is measured by the value
$$\cE_{\pi}(\hat{f}_n)=\mathbb{E}_{\pi}[R(\hfn)-R^*]$$ called \emph{excess
risk of $\hfn$}. In this case $R(\hfn)=\mathbb{P}(\hfn(X)\neq
Y|D_n)$ and $\mathbb{E}_{\pi}$ denotes the expectation w.r.t. $D_n$
when the probability distribution of $(X_i,Y_i)$ is $\pi$ for any
$i=1,\ldots,n$. We say that a classifier $\hf$ learns with the
convergence rate $\phi(n)$, where $(\phi(n))_{n\in\mathbb{N}}$ is a
decreasing sequence, if an absolute constant $C>0$ exists such that
for any integer $n$, $\mathbb{E}_{\pi}[R(\hf) - R^*]\leq C\phi(n)$.

We introduce a loss function on the set of all prediction rules:
$$d_{\pi}(f,g)=\vert R(f)-R(g)\vert.$$ This loss is a
\emph{semi-distance} (it is symmetric, satisfies the triangle
inequality and $d_\pi(f,f)=0$). For all classifiers $\hfn$, it is
linked to the excess risk by
$$\cE_{\pi}(\hfn)=\mathbb{E}_{\pi}[d_\pi(\hfn,f^*)],$$ where the RHS is the risk of $\hfn$
associated to the loss $d_\pi$. In classification we can consider
three estimation problems. The first one is estimation of the Bayes
rule $f^*$, the second one is estimation of the conditional
probability function $\eta$ and the last one is estimation of the
probability $\pi$. Usually, estimation of $\eta$ involves smoothness
assumption on the conditional function $\eta$. However, global
smoothness assumptions on $\eta$ are somehow too restrictive for the
estimation of $f^*$ since the behavior of $\eta$ away from the
decision boundary $\{x\in\cX: \eta(x)=1/2\}$ may have no effect on
the estimation of $f^*$.

In this paper we deal directly with estimation of $f^*$. But, in
this case, the main difficulty of the classification problem is the
dependence on $\pi$ of the loss $d_\pi$ (usually, we use a loss free
from $\pi$, which upper bounds $d_\pi$ to obtain rates of
convergence). Moreover, using the loss $d_\pi$, we don't have the
usual bias/variance trade-off, unlike many other estimation
problems. This is due to the fact that we do not have an
approximation theory in classification for the loss $d_\pi$. This
gap is due to the difficulty that $d_\pi$ depends on $\pi$, thus,
this theory has to be uniform on $\pi$. We need approximation
results of the form:
\begin{equation}\label{approxtheory}
\forall \pi=(P^X,\eta) \in \cP, \forall \epsilon>0, \exists
f_\epsilon\in\cF_\epsilon : d_\pi(f_\epsilon,f^*)\leq \epsilon,
\end{equation}
where $P^X$ is the marginal distribution of $\pi$ on $\cX$,
$f^*=\rm{sign}(2\eta-1)$, $\cP$ is a set of probability measures on
$\cX\times\{-1,1\}$ and the family of classes of prediction rules
$(\cF_\epsilon)_{\epsilon>0}$ is decreasing
($\cF_\epsilon\subset\cF_{\epsilon'}$ if $\epsilon'<\epsilon$) and
$\cF_\epsilon$ is less complex than $\{f^*:\pi\in\cP\}$, in fact we
expect $\cF_\epsilon$ to be parametric. Similar results appear in
density estimation literature, where, for instance, $\cP$ is
replaced by the set of all probability measures with a density with
respect to the Lebegue measure lying in an $L^1-$ball and
$\cF_\epsilon$ is replaced by the set of all functions with a finite
number (depending on $\epsilon$) of coefficients non equal to zero
in the decomposition in the chosen orthogonal basis. But
approximation theory in density estimation does not depend on the
underlying probability measure since the loss functions used there
are generally independent of the underlying statistical problem. In
this paper, we deal directly with the estimation of the Bayes rule
and obtain convergence result w.r.t. the loss $d_\pi$ by using an
approximation approach of the Bayes rules w.r.t. $d_\pi$. Theorems
in Section $7$ of \cite{dgl:96} show that no classifier can learn
with a given convergence rate for arbitrary underlying probability
distribution $\pi$. Thus, assumption on $f^*$ has to be done to
obtain convergence rates. In this paper, assumption on $f^*$ is
close to the one met in density estimation when we assume that the
underlying density belongs to an $L^1-$ball.

Usually, a model (set of measurable functions with values in
$\{-1,1\}$) is considered and we assume that the Bayes rule belongs
to this model. In this case the bias is equal to zero and no bound
on the approximation term is considered. In \cite{blv:03}, question
on the control of the approximation error for a class of models in
the boosting framework is asked. In this paper, it is assumed that
the Bayes rule belongs to the model and nature of distribution
satisfying such condition is explored. Another related work is
\cite{lv:04}, where, under general conditions, it can be guaranteed
that the approximation error converges to zero for some specific
models. In the present paper, bias term is not taken equal to zero
and convergence rates for the approximation error are obtained
depending on the complexity of the considered model (cf. Theorem
\ref{approxtheo}).

We consider the classification problem on $\cX=[0,1]^d$. All the
results can be generalized to a given compact of $\mathbb{R}^d$.
Like in many other works on the classification problem an upper
bound for the loss $d_\pi$ is used. But, in our case we still work
directly with the estimation of $f^*$. For a prediction rule $f$ we
have
\begin{equation}\label{maj}
d_\pi(f_1,f^*)=\mathbb{E}[|2\eta(X)-1|\1_{f_1(X)\neq f^*(X)}]\leq
||f_1-f^*||_{L^1(P^X)}.\end{equation} In order to get a
distribution-free loss function, we assume that the following
assumption holds

\noindent {\it{\bf{(A1)}} The marginal $P^X$ is absolutely
continuous w.r.t. the Lebesgue measure $\lambda_d$ and
 $0<a\leq dP^X(x)/d\lambda_d \leq A <+\infty, \quad
\forall x\in[0,1]^d$.}

This is a technical assumption used for the control of the $P^X$
measure of some subset of $[0,1]^d$. In recent years some
assumptions have been introduced to measure a statistical quality of
classification problems. The behavior of the regression function
$\eta$ near the level $1/2$ is a key point of the classification's
quality (cf. e.g. \cite{tsy:04}). In fact, the closest is $\eta$ to
$1/2$, the more difficult is the classification problem,
nevertheless when we have $\eta\equiv1/2$ the classification is
trivial since all prediction rules are Bayes rules. Here, we measure
the quality of the classification problem thanks to the following
assumption introduced by \cite{mn:03}:

 \noindent {\bf Strong Margin Assumption
(SMA):} There exists an absolute constant $0<h\leq1$ such
that:$$\mathbb{P}\left(\vert 2\eta(X)-1\vert>h \right)=1.$$ Under
assumptions (A1) and (SMA) we have
$$ah||f_1-f^*||_{L_1(\lambda_d)}\leq d_\pi(f_1,f^*)\leq
A||f_1-f^*||_{L_1(\lambda_d)}.$$Thus, estimation of $f^*$ w.r.t. the
loss $d_\pi$ is the same as estimation w.r.t. $L_1(\lambda_d)-$norm,
where $\lambda_d$ is the Lebesgue measure on $[0,1]^d$.

The paper is organized as follows. In the next section we propose a
representation for functions with values in $\{-1,1\}$ in a
fundamental system of $L^2([0,1]^d)$. The third section is devoted
to approximation and estimation of Bayes rules having a sparse
representation in this system. In the fourth section we discuss
about this approach. Proofs are given in the last section.

\par

\section{Classes of Bayes Rules with Sparse Representation}\label{sectionsparseclassofbayesrules}

Theorem \ref{approxtheo} of Subsection
\ref{subsectionapproximationresult} is about the approximation of
the Bayes rules when we assume that $f^*$ belongs to a kind of
"$L^1-$ball" for functions with values in $\{-1,1\}$. The idea is
to develop $f^*$ in a fundamental system of $L^2([0,1]^d,P^X)$
(that is a countable family of functions such that the set of all
finite linear combinations is dense in $L^2([0,1]^d,P^X)$)
inherited from the Haar basis and to control the number of
coefficients non equal to zero. In this paper we only consider the
case where $P^X$ satisfies (A1). We can extend the study to a more
general case by taking another partition of $[0,1]^d$ adapted to
$P^X$.

First we construct such a fundamental system. We consider a
sequence of partitions of $\cX=[0,1]^d$ by setting for any integer
$j$,
$$\cI_{\mathbf{k}}^{(j)}=I_{k_1}^{(j)}\times\ldots\times
I_{k_d}^{(j)},$$ where $\mathbf{k}$ is the multi-index
$$\mathbf{k}=(k_1,\ldots,k_d)\in I_d(j)=\{0,1,\ldots,2^j-1\}^d,$$
and for any integer $j$ and any $k\in\{1,\ldots,2^j-1\},$
$$I_k^{(j)}=\left\{\begin{array}{ll} \left[ \frac{k}{2^j},\frac{k+1}{2^j}\right) & \mbox{if }k=0,\ldots,2^j-2 \\
\left[\frac{2^j-1}{2^j},1 \right] & \mbox{if }k=2^j-1
\end{array}.\right.$$
We consider the family $\cS=\left(
\phi_{\mathbf{k}}^{(j)}:j\in\mathbb{N},\mathbf{k}\in
I_d(j)\right)$ where
$$\phi_{\mathbf{k}}^{(j)}=\1_{\cI_{\mathbf{k}}^{(j)}}, \quad \forall
j\in\mathbb{N}, \mathbf{k}\in I_d(j),$$ where $\1_A$ denotes the
indicator of a set $A$. Set $\cS$ is a fundamental system of
$L^2([0,1]^d,P^X)$. This is the class of indicators of the dyadic
sets of $[0,1]^d$.

We consider the class of functions $f$ defined $P^X-a.s.$ from
$[0,1]^d$ to $\{-1,1\}$ which can be written in this system by
$$
f=\sum_{j=0}^{+\infty}\sum_{\mathbf{k}\in
I_d(j)}a_{\mathbf{k}}^{(j)}\phi_{\mathbf{k}}^{(j)}, P^X-a.s., \mbox{
where } a_{\mathbf{k}}^{(j)}\in\{-1,0,1\},$$ where, for any point
$x\in[0,1]^d$, the right hand side applied in $x$ is a finite sum.
Denote this class by $\cF^{(d)}$. In what follows, we use the
vocabulary appearing in the wavelet literature. The index "$j$" of
$a_{\mathbf{k}}^{(j)}$ and $\phi_{\mathbf{k}}^{(j)}$ is called
"level of frequency". Since $\cS$ is not an orthogonal basis of
$L^2([0,1]^d,P^X)$, the expansion of $f$ w.r.t. this system is not
unique. Therefore, to avoid any ambiguity, we define an unique
writing for any mapping $f$ in $\cF^{(d)}$ by taking
$a_{\mathbf{k}}^{(j)}\in\{-1,1\}$ with preferences for low
frequencies when it is possible.  Roughly speaking, for
$f\in\cF^{(d)}$, denoted by
$f=\sum_{j=0}^{+\infty}\sum_{\mathbf{k}\in
I_d(j)}a_{\mathbf{k}}^{(j)}\phi_{\mathbf{k}}^{(j)}, P^X-a.s.$ where
$a_{\mathbf{k}}^{(j)}\in\{-1,0,1\}$, it means that, we construct
$A_{\mathbf{k}}^{(j)}\in\{-1,0,1\},j\in\mathbb{N},\mathbf{k}\in
I_d(j)$, such that, if there exists $J\in\mathbb{N}$ and
$\mathbf{k}\in I_d(J)$ such that for all $\mathbf{k}'\in I_d(J+1)$
satisfying $\phi_{\mathbf{k}}^{(J)}\phi_{\mathbf{k}'}^{(J+1)}\neq0$
we have $a_{\mathbf{k}'}^{(J+1)}=1$, then we take
$A_{\mathbf{k}'}^{(J)}=1$ and the $2^d$ other coefficients of higher
frequency $A_{\mathbf{k}'}^{(J+1)}=0$ instead of having these $2^d$
coefficients equal to $1$, and the same convention holds for $-1$.
Moreover if we have $A_{\mathbf{k}}^{(J_0)}\neq0$ then
$A_{\mathbf{k}'}^{(J)}=0$ for all $J>J_0$ and $\mathbf{k}'\in
I_d(J)$ satisfying
$\phi_{\mathbf{k}}^{(J_0)}\phi_{\mathbf{k}'}^{(J)}\neq0$. We can
describe a mapping $f\in\cF^{(d)}$ satisfying this convention by
using a tree. Each knot corresponds to a coefficient
$A_{\mathbf{k}}^{(J)}$. The root is $A_{0,\ldots,0}^{(0)}$. If a
knot, describing the coefficient $A_{\mathbf{k}}^{(J)}$, equals to
$1$ or $-1$ then it has no branches, otherwise it has $2^d$
branches, corresponding to the $2^d$ coefficients at the following
frequency, describing the coefficients $A_{\mathbf{k}'}^{(J+1)}$ for
$\mathbf{k}'$ satisfying
$\phi_{\mathbf{k}}^{(J)}\phi_{\mathbf{k}'}^{(J+1)}\neq0$. At the end
all the leaves of the tree equals to $1$ or $-1$, and the depth of a
leaf is the frequency of the coefficient associated. The writing
convention says that a knot can not have all his leaves equal to $1$
together (or $-1$). In this case we write this mapping by putting a
$1$ at the knot (or $-1$). In what follows we say that a function
$f\in\cF^{(d)}$ satisfies the writing convention (W) when $f$ is
written in $\cS$ using the writing convention describes in this
paragraph. Remark that this writing convention is not an assumption
on the function since we can write all $f\in\cF$ using this
convention. Representation of the Bayes rules using Dyadic decision
trees has been explored by \cite{ns:04}.

Is it possible to write every measurable functions from $[0,1]^d$
to $\{-1,1\}$ in the fundamental system $\cS$ using coefficients
with values in $\{-1,0,1\}$? Since the family of set
$(\mathring{\cI}_{\mathbf{k}}^{(j)}:j\in\mathbb{N},\mathbf{k}\in
I_d(j))$, where $\mathring{A}$ denotes the interior of $A$, is a
basis of open subsets of $[0,1]^d$, this question is equivalent to
this one: "Take $A$ a Borel of $[0,1]^d$, is it possible to find
an open subset $\cO$ of $[0,1]^d$ such that the symmetrical
difference between $A$ and $\cO$ has a Lebesgue measure $0$?"
Unfortunately, the answer to this last question is negative. There
exists $F\subset[0,1]^d$ a Borel,  closed, with an empty interior
and a positive Lebesgue measure $\lambda_{d}(F)>0$. For example,
in the one dimension case, the following algorithm yields such a
set. Take $(l_k)_{k\geq1}$ a sequence of numbers defined by
$l_k=1/2-1/(k+1)^2$ for any integer $k$. Denote by $F_0$ the
interval $[0,1]$ and construct a sequence of closed sets
$(F_k)_{k\geq0}$ like in the following picture.

{\centerline{\unitlength 1cm
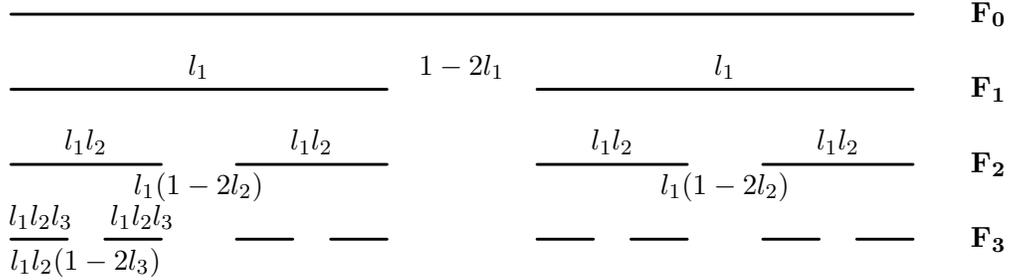
\begin{picture}(16,6.95)(-7.5,-2.45)
\Thicklines 
\path(-6,2)(6,2)(-6,2)
\put(7,2){\makebox(0,0){\small $\mathbf{F_0}$}}
\path(-6,1)(-1,1)(-6,1)
\path(1,1)(6,1)(1,1)
\put(7,1){\makebox(0,0){\small $\mathbf{F_1}$}}
\path(-6,-0)(-4,-0)(-6,-0)
\path(-3,-0)(-1,-0)(-3,-0)
\path(1,0)(3,0)(1,0)
\path(4,0)(6,0)(4,0)
\put(7,0){\makebox(0,0){\small $\mathbf{F_2}$}}
\path(-6,-1)(-5.25,-1)(-6,-1)
\path(-4.75,-1)(-4,-1)(-4.75,-1)
\path(-3,-1)(-2.25,-1)(-3,-1)
\path(-1.75,-1)(-1,-1)(-1.75,-1)
\path(1,-1)(1.75,-1)(1,-1)
\path(2.25,-1)(3,-1)(2.25,-1)
\path(4,-1)(4.75,-1)(4,-1)
\path(5.25,-1)(6,-1)(5.25,-1)
\put(7,-1){\makebox(0,0){\small $\mathbf{F_3}$}}
\put(-3.5,1.3){\makebox(0,0){\small $l_1$}}
\put(3.5,1.3){\makebox(0,0){\small $l_1$}}
\put(-5,0.3){\makebox(0,0){\small $l_1 l_2$}}
\put(0,1.3){\makebox(0,0){\small $1-2l_1$}}
\put(-3.5,-0.3){\makebox(0,0){\small $l_1(1-2l_2)$}}
\put(-2,0.3){\makebox(0,0){\small $l_1l_2$}}
\put(2,0.3){\makebox(0,0){\small $l_1 l_2$}}
\put(3.5,-0.3){\makebox(0,0){\small $l_1(1-2l_2)$}}
\put(5,0.3){\makebox(0,0){\small $l_1l_2$}}
\put(-5.6,-0.7){\makebox(0,0){\small $l_1l_2l_3$}}
\put(-5,-1.3){\makebox(0,0){\small $l_1l_2(1-2l_3)$}}
\put(-4.25,-0.7){\makebox(0,0){\small $l_1l_2l_3$}}
\end{picture}
}}

 It is easy to check that
$F=\cap_{k\geq0}F_k$ is closed, with an empty interior and a
positive Lebesgue measure. For the $d$-dimensional case, the set
$F\times[0,1]^{d-1}$ satisfies the required assumptions. Thus, take
$F$ such a set and $\cO$ an open subset of $[0,1]^d$. If
$\cO\subseteq F$ then $\cO=\emptyset$ because
$\mathring{F}=\emptyset$ and
$\lambda_{d}(F\Delta\cO)=\lambda_{d}(F)>0$. If $\cO\not\subseteq F$
then $\cO\cap F^{c}$ is an open subset of $[0,1]$ none empty, so
$\lambda_{d}(\cO\Delta F)\geq\lambda_d(\cO\cap F^{c})>0$. Thus,
every measurable functions from $[0,1]^d$ to $\{-1,1\}$ can not be
written in $\cS$ using only coefficients with values in
$\{-1,0,1\}$. Nevertheless, the Lebesgue measure satisfies the
property of regularity, which says that for any Borel $B\in[0,1]^d$
and any $\epsilon>0$, there exists a compact subset $K$ and an open
subset $\cO$ such that $K\subseteq A \subseteq \cO$ and
$\lambda_{d}(\cO-K)\leq\epsilon$. Hence, one can easily check that
for any measurable function $f$ from $[0,1]^d$ to $\{-1,1\}$ and any
$\epsilon>0$, there exists a function $g\in\cF^{(d)}$ such that
$\lambda_d(\{x\in[0,1]^d:f(x)\neq g(x)\})\leq \epsilon$. Thus,
$\cF^{(d)}$ is dense in $L^{2}(\lambda_d)$ intersected with the set
of all measurable functions from $[0,1]^d$ to $\{-1,1\}$. Now, we
exhibit some usual prediction rules which belong to $\cF^{(d)}$.

\begin{defi}
  Let $A$ be a Borel subset of $[0,1]^d$.
  We say that $A$ is
  {\bf almost everywhere open} if there exists an open subset $\cO$ of
  $[0,1]^d$ such that $\lambda_d(A\Delta \cO)=0$, where $\lambda_d$ is
  the Lebesgue measure on $[0,1]^d$ and $A\Delta\cO$
  is the symmetrical difference.
\end{defi}

\begin{theo}\label{theoopenas}
  Let $\eta$ be a function from $[0,1]^d$ to $[0,1]$. We consider
  \begin{equation*}
    f_\eta(x)=\left\{
        \begin{array}{cc}
            1 & \mbox{ if } \eta(x)\geq1/2\\
            -1 & \mbox{ otherwise.}
        \end{array}
    \right.
  \end{equation*}
  We assume that $\{\eta\geq1/2\}$ and $\{\eta<1/2\}$ are almost
  everywhere open. Thus, there exists $g\in\cF$ such that for
  $\lambda_d$-almost every
   $x\in[0,1]^d, \quad g=f_\eta, \lambda_d-a.s.$. For
   instance, if $\lambda_d(\partial \{\eta=1/2\})=0$ and, either $\eta$ is $\lambda_d$-almost everywhere
   continuous (it means that there exists an open subset of $[0,1]^d$ with a
   Lebesgue measure equals to $1$ such that
   $\eta$ is continuous on this open subset) or if $\eta$ is $\lambda_d-$almost everywhere equal to a continuous
   function, then $f_\eta\in\cF^{(d)}$.
\end{theo}

Now, we define a model for the Bayes rule by taking a subset of
$\cF^{(d)}$. For all functions $w$ defined on $\mathbb{N}$ and
with values in $\mathbb{R}^+$, we consider  $\cF_w^{(d)}$, the
model for Bayes rules, made of all prediction rules $f$ which can
be written, using the previous writing convention (W), by
$$
f=\sum_{j=0}^{+\infty}\sum_{\mathbf{k}\in
I_d(j)}a_{\mathbf{k}}^{(j)}\phi_{\mathbf{k}}^{(j)},$$ where
$a_{\mathbf{k}}^{(j)}\in\{-1,0,1\}$  and
$$\rm{card}\left\{ \mathbf{k}\in
I_d(j):a_{\mathbf{k}}^{(j)}\neq 0\right\}\leq w(j), \quad \forall
j\in\mathbb{N}.$$ The class $\cF_w^{(d)}$ depends on the choice of
the function $w$. If $w$ is too small then the class $\cF_w^{(d)}$
is not very rich, that is the subject of the following Proposition
\ref{p1}. If $w$ is too large then $\cF_w^{(d)}$ would be too
complex for a good estimation of $f^*\in\cF_w^{(d)}$, that is why
we introduce Definition \ref{elipsoide} in what follows.

\begin{pro}\label{p1} Let $w$ be a mapping from $\mathbb{N}$ to
$\mathbb{R}^+$ such that $w(0)\geq 1$. The two following
assertions are equivalent:
\begin{enumerate}[(i)]
\item $\cF_w^{(d)}\neq \{\1_{[0,1]^d}\}$. \item $
\sum_{j=1}^{+\infty} 2^{-dj}\lfloor w(j) \rfloor\geq1.$
\end{enumerate}
\end{pro}

And if $w$ is too large then the approximation by a parametric
model will be impossible, that is why we give a particular look on
the class of function introduced in the following Definition
\ref{elipsoide}.

\begin{defi}\label{elipsoide} Let $w$ be a mapping from $\mathbb{N}$ to
$\mathbb{R}^+$. If $w$ satisfies
\begin{equation}\label{L1Ball}
\sum_{j=0}^{+\infty} \frac{\lfloor w(j)
\rfloor}{2^{dj}}<+\infty,\end{equation} then we say that
$\cF_w^{(d)}$ is a {\bf $\mathbf{L^1-}$ball of prediction rules}.
\end{defi}

\begin{rem}
We say that $\cF_w^{(d)}$ is a "$L^1-$ball" for a function $w$
satisfying (\ref{L1Ball}), because , the sequence $(\lfloor w(j)
\rfloor)_{j\in\mathbb{N}}$ belongs to a $L^1-$ball of
$\mathbb{N}^\mathbb{N}$, with radius $(2^{dj})_{j\in\mathbb{N}}$.
Moreover, definition \ref{elipsoide} can be link to the definition
of a $L^1-$ball for real valued functions, since we have a kind of
base, given by $\cS$, and we have a control on coefficients which
increases with the frequency. Control on coefficients, given by
(\ref{L1Ball}), is close to the one for coefficients of a real
valued function in $L^1-$ball since it deals with the quality of
approximation of the class $\cF_w^{(d)}$ by a parametric model.
\end{rem}

\begin{rem}
  A $L^1-$ball of prediction rules is made of "sparse"
  prediction
  rules. In fact, for $f\in\cF_w^{(d)}$, the repartition of
  coefficients non equal to zero in the decomposition of $f$ at a
  given frequency becomes sparse as the frequency grows. That is
  the reason why $\cF_w^{(d)}$ can be called a {\bf{sparse class of
  prediction rules}}. For exemple, if $(\lfloor w(j)
\rfloor/2^{dj})_{j\geq1}$ decreases and (\ref{L1Ball}) holds then
number of coefficients non equal to $0$ at the frequency $j$ is
smaller than $j^{-1}$ per cent of the maximal number of
coefficients (that is $2^{dj}$).
\end{rem}

\begin{rem} If we assume that $P^X$ is known then we can
work with any measurable space $\cX$ endowed with a Lebesgue measure
$\lambda$, while assuming that $P^X<<\lambda$. In this case, we take
$ \left(\cI_{\mathbf{k}}^{(j)}:j\in\mathbb{N}, \mathbf{k}\in
I_d(j)\right)$, such that for any $j\in\mathbb{N}$,
$\left(I_{\mathbf{k}}^{(j)}:\mathbf{k}\in I_d(j)\right)$ is a
partition of $\cX$ adapted to the previous one $
\left(I_{\mathbf{k}}^{(j-1)}: \mathbf{k}\in I_d(j-1)\right)$ and
satisfying $P^X(I_{\mathbf{k}}^{(j)})=2^{-jd}$. All the results
below can be obtained in this framework.
\end{rem}

Now, examples of functions satisfying (\ref{L1Ball}) are given.
Classes $\cF_w^{(d)}$ associated to these functions are used in what
follows to define statistical models. As an introduction we define
the minimal infinite class of prediction rules, by $\cF_0^{(d)}$
which is the class $\cF_w^{(d)}$ for $w=w_0^{(d)}$ where
$w_0^{(d)}(0)=1$ and $w_0^{(d)}(j)=2^d-1$, for all $j\geq1$. To
understand why this class is important we introduce a notion of
local oscillation of a prediction rule. This concept defines a kind
of "regularity" for functions with values in $\{-1,1\}$.

\begin{defi} Let $f$ be a prediction rule from $[0,1]^d$ to
$\{-1,1\}$ in $\cF^{(d)}$. We consider the writing of $f$ in the
fundamental system introduce in Section
\ref{subsectionapproximationresult} with writing convention (W):
$$f=\sum_{j=0}^{+\infty}\sum_{\mathbf{k}\in I_d(j)}a_{\mathbf{k}}^{(j)}
\phi_{\mathbf{k}}^{(j)}, P^X-a.s..$$ Let $J\in\mathbb{N}$ and
$\mathbf{k}\in I_d(j)$. We say that $I_{\mathbf{k}}^{(J)}$ is a {\bf
low oscillating block} of $f$ when  $f$ has exactly $2^d-1$
coefficients, in this block, non equal to zero at each level of
frequencies greater than $J+1$. In this case we say that {\bf $f$
has a low oscillating block of frequency $J$}.
\end{defi}

Remark that, if $f$ has an oscillating block of frequency $J$,
then $f$ has an oscillating block of frequency $J'$, for all
$J'\geq J$. The function class $\cF_0^{(d)}$ is made of all
prediction rules with one oscillate block at level $1$ and of the
indicator function $\1_{[0,1]^d}$. If we have
$w(j_0)<w_0^{(d)}(j_0)$ for one $j_0\geq1$ and $w(j)=w_0^{(d)}(j)$
for $j\neq j_0$ then the associated class $\cF_w^{(d)}$ contains
only the indicator function $\1_{[0,1]^d}$, that is the reason why
we say that $\cF_0^{(d)}$ is "minimal".

Nevertheless,  the following proposition shows that $\cF_0^{(d)}$
is a rich class of prediction rules from a combinatorial point of
view. We recall some quantities which measure a combinatorial
richness of a class of prediction rules. For any class $\cF$ of
prediction rules from $\cX$ to $\{-1,1\}$, we consider
$$N(\cF,(x_1,\ldots,x_m))=card\left(\{(f(x_1),\ldots,f(x_m)):f\in\cF\}
\right)$$ where $x_1,\ldots,x_m\in\cX$ and $m\in\mathbb{N}$,
$$S(\cF,m)=\max\left(
N(\cF,(x_1,\ldots,x_m)):x_1,\ldots,x_m\in\cX\right)$$ and the
$VC$-dimension of $\cF$ is
$$VC(\cF)=\min\left(m\in\mathbb{N}:S(\cF,m)\neq 2^m \right).$$
Consider
$x_j=\left(\frac{2^j+1}{2^{j+1}},\frac{1}{2^{j+1}},\ldots,\frac{1}{2^{j+1}}
\right),$ for any $j\in\mathbb{N}$. Thus, for any integer $m$, we
have $N(\cF_0^{(d)},(x_1,\ldots,x_m))=2^m$. Hence, the following
proposition holds.
\begin{pro}\label{l2}
The class of prediction rules $\cF_{0}^{(d)}$ has an infinite
$VC$-dimension.
\end{pro}

Thus every class $\cF_w^{(d)}$ such that $w\geq w_0^{(d)}$ has an
infinite $VC$-dimension (since $w\leq w'\Rightarrow
\cF_w^{(d)}\subseteq\cF_{w'}^{(d)}$), which is the case for the
following classes.

Now, we introduce some examples of $L^1-$ball of Bayes rules. We
denote by $\cF_K^{(d)}$, for a $K\in\mathbb{N}^*$, the class
$\cF_w^{(d)}$ of prediction rules where $w$ is equal to the
function
$$ w_K^{(d)}(j)=\left\{\begin{array}{ll} 2^{dj} & \mbox{if }  j\leq K, \\
2^{dK} & \mbox{otherwise. } \end{array}\right.$$ This class is
called the {\bf truncated class of level K}.

We consider {\bf exponential classes}. These sets of prediction
rules are denoted by $\cF_{\alpha}^{(d)}$, where $0<\alpha <1$,
and are equal to $\cF_w^{(d)}$ when $w=w_{\alpha}^{(d)}$ and
$$w_{\alpha}^{(d)}(j)=
\left\{\begin{array}{ll} 2^{dj} & \mbox{if }  j\leq N^{(d)}(\alpha) \\
2^{d\alpha j} & \mbox{otherwise } \end{array}\right. ,$$where
$N^{(d)}(\alpha)=\inf\left(N\in\mathbb{N}:2^{d\alpha N}\geq 2^d-1
\right)$, that is for $N^{(d)}(\alpha)=\lceil
\log(2^d-1)/(d\alpha\log 2)\rceil$.

\begin{rem}
  For the one-dimensional case, an other point of view is to consider $f^*\in L^2([0,1])$ and
  to develop $f^*$ in an orthogonal basis of $L^2([0,1])$.
  Namely,
  $$f^*=\sum_{j\in\mathbb{N}}\sum_{k=0}^{2^j-1}a_k^{(j)}\psi_k^{(j)},$$where
  $a_k^{(j)}=\int_0^1 f^*(x)\psi_k^{(j)}(x)dx$ for any
  $j\in\mathbb{N}$ and $k=0,\ldots,2^j-1$. For the control of the
  bias term we assume that the family of
  coefficients $(a_k^{(j)},j\in\mathbb{N},k=0,\ldots,2^j-1)$
  belongs to a $L^1-$ball. But this point of view leads to analysis and estimation issues.
  First problem: Which functions with values in $\{-1,1\}$ have wavelet coefficients
  in a $L^1-$ball and which wavelet basis is more adapted to our problem (maybe the Haar basis)?
   Second problem: Which kind of estimators could
  be used for the estimation of these coefficients? As we can see,
  the main problem is that there is no approximation theory for
  functions with values in $\{-1,1\}$. We do not know how to
  approach, in $L^2([0,1])$, measurable functions with values in
  $\{-1,1\}$ by "parametric" functions with values in $\{-1,1\}$.
  Methods developed in this paper may be seen as a first step in
  this field. We can generalize this approach to functions with
  values in $\mathbb{Z}$. Remark that when functions take values
  in $\mathbb{R}$, that is for the regression problem, usual
  approximation theory is used to obtain a control on the bias
  term.
\end{rem}

\begin{rem}
  Other sets of prediction rules are described by the classes
  $\cF_w^{(d)}$ where $w$ is from $\mathbb{N}$ to $\mathbb{R}^+$
  and satisfies $$\sum_{j\geq1}a_j\frac{\lfloor w(j)
\rfloor}{2^{dj}}\leq
  L,$$where $(a_j)_{j\geq1}$ is an increasing sequence of positive
  numbers.
\end{rem}

\section{Rates of Convergence over $\cF^{(d)}_w$ under
(SMA)}\label{sectionrateofconvergence}
\subsection{Approximation Result}\label{subsectionapproximationresult}

Let $w$ be a function from $\mathbb{N}$ to $\mathbb{R}^+$ and
$A>1$, we denote by $\cP_{w,A}$ the set of all probability
measures $\pi$ on $[0,1]^d\times\{-1,1\}$ such that the Bayes
rules $f^*$, associated to $\pi$, belongs to $\cF_w^{(d)}$ and the
marginal of $\pi$ on $[0,1]^d$ is absolutely continuous and one
version of its Lebesgue density is upper bounded by $A$. The
following Theorem can be seen as an approximation Theorem for the
Bayes rules w.r.t. the loss $d_\pi$ uniformly in
$\pi\in\cP_{w,A}$.

\begin{theo}[Approximation Theorem]\label{approxtheo}
Let $\cF_w^{(d)}$ be a $L^1-$ball of prediction rules. We have:
$$\forall \epsilon>0, \exists J_\epsilon\in\mathbb{N}:\forall \pi\in\cP_{w,A},
 \exists
 f_\epsilon=\sum_{\mathbf{k}\in I_d(J_\epsilon)}B_{\mathbf{k}}^{(J_\epsilon)}
\phi_{\mathbf{k}}^{(J_\epsilon)}$$  where
$B_{\mathbf{k}}^{(J_\epsilon)}\in\{-1,1\}$ and $$
d_{\pi}(f^*,f_\epsilon)\leq \epsilon,$$where $f^*$ is the Bayes rule
associated to $\pi$. For example, $J_\epsilon$ can be the smallest
integer $J$ satisfying $\sum_{j=J+1}^{+\infty}2^{-dj}\lfloor
w(j)\rfloor<\epsilon/A.$
\end{theo}

\begin{rem} No assumption on the quality of the classification
problem, like an assumption on the margin, is needed to state
Theorem \ref{approxtheo}. Only assumption on the "number
 of oscillations" of $f^*$ is used. Theorem
 \ref{approxtheo} deals with approximation of functions in the $L^1-$ball $\cF_w^{(d)}$ by functions with values in
 $\{-1,1\}$ and no
 estimation issues are met.
\end{rem}

\begin{rem}
Theorem \ref{approxtheo} is the first step to prove an estimation
theorem using a trade-off between a bias term and a variance term.
We write
$$\cE_{\pi}(\hfn)=\mathbb{E}_{\pi}\left[d_\pi(\hfn,f^*)
\right]\leq
\mathbb{E}_{\pi}\left[d_\pi(\hfn,f_\epsilon)\right]+d_\pi(f_\epsilon,f^*).$$
Since $f_\epsilon$ belongs to a parametric model we expect to have
a control of the variance term,
$\mathbb{E}_{\pi}\left[d_\pi(\hfn,f_\epsilon)\right]$, depending
on the dimension of the parametric model which is linked to the
quality of the approximation in the bias term.
\end{rem}

\begin{rem}
Since
$d_\pi(f^*,f_\epsilon)=\mathbb{E}\left[|2\eta(X)-1|\1_{f^*(X)\neq
f_\epsilon(X)}\right]$, the closest to $1/2$  $\eta$ is, the
smallest the bias is. Especially, we have a bias equal to zero
when $\eta=1/2$ (in this case any prediction rule is a Bayes
rules). Thus, more difficult the problem of estimation is (that is
for underlying probability measure $\pi=(P^X,\eta)$ with $\eta$
close to $1/2$), the smallest the bias is. This behavior does not
appear clearly in density estimation.
\end{rem}

\par

\subsection{Estimation Result} We consider the following class of
estimators indexed by the frequency rank $J\in\mathbb{N}$:
\begin{equation}\label{esti}\hfn^{(J)}=\sum_{\mathbf{k}\in I_d(J)}\hat{A}_{\mathbf{k}}^{(J)}
\phi_{\mathbf{k}}^{(J)},\end{equation} where coefficients are
defined by$$ \hat{A}_{\mathbf{k}}^{(J)}=\left\{\begin{array}{ll} 1 &
\mbox{if } \exists X_i\in I_{\mathbf{k}}^{(J)} \mbox{ and
}\rm{card}\left\{i:\begin{array}{l} X_i\in
I_{\mathbf{k}}^{(J)},\\Y_i=1\end{array} \right\}>\rm{card}\left\{i:
\begin{array}{l}X_i\in
I_{\mathbf{k}}^{(J)},\\Y_i=-1\end{array} \right\} \\
-1 & \mbox{otherwise} \end{array},\right.$$

To obtain a good control of the variance term, we need to assure a
good quality of the estimation problem. Therefore, estimation
results are obtained in Theorem \ref{estitheo} under (SMA)
assumption. In recent years we have understood that (SMA)
assumption can lead to fast rates but is not enough to assure any
rate of convergence (cf. corolary \ref{theodgl} at the end of
section \ref{subsectionoptimality}), thus we have to define a
model for $\eta$ or $f^*$, here we use a $L^1-$ball of prediction
rules as a model for $f^*$.

\begin{theo}[estimation Theorem]\label{estitheo} Let $\cF_w^{(d)}$ be a
$L^1-$ball of prediction rules. Let $\pi$ be a probability measure
on $[0,1]^d\times\{-1,1\}$ satisfying  assumptions (A1) and (SMA),
and such that the Bayes rule, associated to $^pi$, belongs to
$\cF_w^{(d)}$. The excess risk of the classifier
$\hfn^{(J_\epsilon)}$ satisfies for any positive number $\epsilon$,
$$\cE_{\pi}(\hfn^{(J_\epsilon)})=\mathbb{E}_{\pi}\left[d_\pi(\hfn^{(J_\epsilon)},f^*) \right]\leq
(1+A)\epsilon+\exp\left(-na(1-\exp(-h^2/2))2^{-dJ_\epsilon}
\right),$$ where $J_\epsilon$ is the smallest integer satisfying
$\sum_{j=J_\epsilon+1}^{+\infty}2^{-dj}\lfloor
w(j)\rfloor<\epsilon/A$. Parameters $a,A$ appear in Assumption
(A1) and $h$ is used in (SMA).
\end{theo}

\begin{rem}
The upper bound can be split in the bias term: $\epsilon$ and the
variance term:
$A\epsilon+\exp\left(-na(1-\exp(-h^2/2))2^{-dJ_\epsilon} \right)$.
Remark that a bias term appears in the variance term.
\end{rem}

\par

\subsection{Optimality}\label{subsectionoptimality}
This section is devoted to the optimality, in a minimax sense, of
estimation in classification models such that $f^*\in\cF_w^{(d)}$.
Let $0<h<1$, $0<a\leq1\leq A <+\infty$ and $w$ a mapping from
$\mathbb{N}$ to $\mathbb{R}^+$. we denote by $\cP_{w,h,a,A}$ the
set of all probability measures $\pi=(P^X,\eta)$ on $[0,1]^d\times
\{-1,1\}$ such that
\begin{enumerate} \item The marginal $P^X$ satisfies (A1). \item
The Assumption (SMA) is satisfied. \item The Bayes rule $f^*$,
associated to $\pi$, belongs to $\cF_w^{(d)}$.\end{enumerate} We
use the version of Lemma of Assouad in the appendix of
\cite{lec3:06} to lower bound the minimax risk on $\cP_{w,h,a,A}$.
From Theorem~\ref{estitheo} and Theorem~\ref{optimal}, we can
deduce the optimality (up to a logarithm term) of the estimator
$\hfn^{(J_n)}$ where the rank $J_n$ is obtained by an optimal
trade-off between the bias term and the variance term.

\begin{theo}\label{optimal} Let $w$ be a function from
$\mathbb{N}$ to $\mathbb{R}^+$ such that \begin{enumerate}[(i)]
\item $\lfloor w(0)\rfloor\geq 1$ and $\forall j\geq1, \quad
\lfloor w(j)\rfloor\geq 2^d-1$ \item $\forall j\geq 1, \quad
\lfloor w(j-1)\rfloor\geq 2^{-d}\lfloor w(j)\rfloor$.
\end{enumerate} We have for all $n\in\mathbb{N}$,
$$ \inf_{\hat{f}_n}\sup_{\pi\in\cP_{w,h,a,A}}
\cE_{\pi}(\hat{f}_n)\geq C_0 n^{-1}\left(\lfloor w\left(\lfloor
\log n/(d\log2)\rfloor+1\right)\rfloor-(2^d-1)\right),$$ and if
$\lfloor w(j)\rfloor\geq 2^d,\quad \forall j\geq1$ then $
\inf_{\hat{f}_n}\sup_{\pi\in\cP_{w,h,a,A}}
\cE_{\pi}(\hat{f}_n)\geq C_0 n^{-1}$ where $C_0=(h/8)\exp\left(
-(1-\sqrt{1-h^2})\right)$.
\end{theo}

\begin{rem} For a function $w$ satisfying assumptions of Theorem \ref{optimal} and under (SMA),
we can not expect a convergence rate faster than $1/n$, which is
the usual lower bound for the classification problem under (SMA).
\end{rem}

From the previous Theorem we obtain immediately Theorem 7.1 of
\cite{dgl:96}. We denote by $\cP_{1}$ the class of all probability
measures on $[0,1]^d\times\{-1,1\}$ such that the marginal
distribution $P^X$ is $\lambda_d$ (the Lebesgue probability
distribution on $[0,1]^d$) and (SMA) is satisfied with the margin
$h=1$. The case "$h=1$" is equivalent to $R^*=0$. That is for a
perfect classification problem, where $Y$ is an exact function of
$X$ given by $Y=f^*(X)=\eta(X)$.

\begin{coro}\label{theodgl} For any integer $n$ we have
$$\inf_{\hat{f}_n}\sup_{\pi\in\cP_{1}}
\cE(\hat{f}_n)\geq \frac{1}{8e}.$$
\end{coro}
It means that no classifier can achieve a rate of convergence in
the classification models $\cP_{1}$, even if these classification
problems are all very good ($Y$ is given by $f^*(X)$ without any
noise and there are no spot of low probability).

\par

\subsection{Rates of Convergence for Different Classes of
Prediction Rules} In this section we apply  results stated in
Theorem \ref{estitheo} and Theorem~\ref{optimal} to different
$L^1-$ball classes $\cF_w^{(d)}$ introduced at the end of Section
\ref{sectionsparseclassofbayesrules}. We give rates of convergence
and lower bounds for these models. Using notations introduced in
Section \ref{sectionsparseclassofbayesrules} and subsection
\ref{subsectionoptimality}, we consider the following models. For
$w=w_{K}^{(d)}$ denote by $\cP_{K}^{(d)}$ the set
$\cP_{w_{K}^{(d)},h,a,A}$ of probability measures on
$[0,1]^d\times\{-1,1\}$ and $\cP_{\alpha}^{(d)}$ for
$w=w_\alpha^{(d)}$.

\begin{theo}\label{t2}
For the truncated class $\cF_K^{(d)}$, we have
$$\sup_{\pi\in\cP_{K}^{(d)}}\cE_{\pi}(\hfn^{(J_n)})\leq
C_{K,h,a,A}\frac{\log n}{n},$$ where $C_{K,h,a,A}>0$ is depending
only on $K,h,a,A$ and for the lower bound, there exists
$C_{0,K,h,a,A}>0$ depending only on $K,h,a,A$ such that, for all
$n\in\mathbb{N}$,
$$ \inf_{\hat{f}_n}\sup_{\pi\in\cP_{K}^{(d)}}
\cE_{\pi}(\hat{f}_n)\geq C_{0,K,h,a,A} n^{-1}.$$

For the exponential class $\cF_{\alpha}^{(d)}$ where $0<\alpha<1$,
we have for any integer $n$
$$
\sup_{\pi\in\cP_{\alpha}^{(d)}}\cE_{\pi}(\hfn^{(J_n)})\leq
C'_{\alpha,h,a,A}\left(\frac{\log n}{n}\right)^{1-\alpha},$$ where
$C'_{\alpha,h,a,A}>0$ and for the lower bound, there exists
$C'_{0,\alpha,h,a,A}>0$ depending only on $\alpha,h,a,A$ such
that, for all $n\in\mathbb{N}$,
$$ \inf_{\hat{f}_n}\sup_{\pi\in\cP_{\alpha}^{(d)}}
\cE_{\pi}(\hat{f}_n)\geq C'_{0,\alpha,h,a,A} n^{-1+\alpha}.$$

In both classes, order of  $J_n$ is $\lceil \log \left(an/(2^d
\log n)\right)/(d\log 2)\rceil$, up to a multiplying constant.
\end{theo}

A remarkable point is that the class $\cF_{K}^{(d)}$ has an
infinite VC-dimension (cf. Section
\ref{sectionsparseclassofbayesrules}). Nevertheless, the rate
$\log n/n$ is achieved on this model.

\section{Discussion}

In this section we discuss about representation and estimation of
"simple" prediction rules in our framework. In considering the
classification problem over the square $[0,1]^2$, a classifier has
to be able to approach, for instance, the "simple" Bayes rule
$f^*_{\cC}$ which is equal to $1$ inside $\cC$, where $\cC$ is a
disc of $[0,1]^2$, and $-1$ outside $\cC$. In our framework, two
questions need to be considered:
\begin{itemize}
  \item How is the representation of the simple function $f^*_\cC$
  in our fundamental system, using only coefficients with values
  in $\{-1,0,1\}$ and with the writing convention (W)?
  \item Is the estimate ${\hat{f}}_n^{(J_n)}$, where $J_n=\lceil \log \left(an/(2^d
\log n)\right)/(d\log 2)\rceil$ is the frequency rank appearing in
Theorem \ref{t2}, a good classifier when
  the underlying probability measure has $f^*_\cC$ for Bayes rule?
\end{itemize}

At a first glance, our point of view is not the right way to
estimate $f^*_\cC$. In this regular case (the border is an infinite
differentiable curve), the direct estimation of the border is a
better approach. The main reason is that a $2$-dimensional
estimation problem becomes a $1$-dimensional problem. Such reduction
of dimension makes estimation easier (in passing, our approach is
specifically good in the $1$-dimensional case, since the notion of
border does not exist in this case). Nevertheless, our approach is
applicable for the estimation of such functions (cf. Theorem
\ref{circle}). Actually, direct estimation of the border reduces the
dimension but there is a big waste of observations since
observations far from the border are not used for this estimation
point of view. It may explain why our approach is applicable. Denote
by
$$\cN(A,\epsilon,||.||_{\infty})=\min\left(N:\exists x_1,\ldots,x_N\in\mathbb{R}^2:
A\subseteq\cup_{j=1}^N B_\infty(x_j,\epsilon) \right)$$ the
$\epsilon-$covering number of a subset $A$ of $[0,1]^2$, w.r.t. the
infinity norm of $\mathbb{R}^2$. For example, the circle
$\cC=\{(x,y)\in\mathbb{R}^2:(x-1/2)^2+(y-1/2)^2=(1/4)^2\}$ satisfies
$\cN(\cC,\epsilon,||.||_{\infty})\leq (\pi/4)\epsilon^{-1}.$ For any
set $A$ of $[0,1]^2$, denote by $\partial A$ the border of $A$.

\begin{theo}\label{circle}
Let $A$ be a subset of $[0,1]^2$ such that $\cN(\partial
A,\epsilon,||.||_{\infty})\leq \delta(\epsilon),$ for any
$\epsilon>0$, where $\delta$ is a decreasing function from
$\mathbb{R}^*_+$ with values in $\mathbb{R}^+$ satisfying
$\epsilon^2\delta(\epsilon)\longrightarrow0$ when $\epsilon$ tends
to zero. Consider the prediction rule $f_A=2\1_{A}-1$. For any
$\epsilon>0$, denote by $\epsilon_0$ the greatest positive number
satisfying $\delta(\epsilon_0)\epsilon_0^2\leq \epsilon$. There
exists a prediction rule constructed in the fundamental system $\cS$
at the frequency rank $J_{\epsilon_0}$ with coefficients in
$\{-1,1\}$ denoted by
$$f_{\epsilon_0}=\sum_{\mathbf{k}\in I_2(J_{\epsilon_0})}
a_{\mathbf{k}}^{(J_{\epsilon_0})}\phi^{(J_{\epsilon_0})}_{\mathbf{k}},$$
with $J_{\epsilon_0}=\lfloor\log(1/\epsilon_0)/\log2\rfloor$ such
that
$$||f_{\epsilon_0}-f_A||_{L^1(\lambda_2)}\leq 36\epsilon.$$
\end{theo}

For instance, there exists a function $f_n$, written in the
fundamental system $\cS$ at the frequency level
$J_n=\lfloor\log(4n/(\pi\log n))/\log2\rfloor$, which approaches the
prediction rule $f_\cC$ with a $L^1(\lambda_2)$ error upper bounded
by $36(\log n)/n$. This frequency level is, up to a multiplying
constant, the same one appearing in Theorem \ref{t2}. In a more
general way, any prediction rule with a border having a finite
perimeter (for instance polygons) is approached by a function
written in the fundamental system at the same frequency rank $J_n$
and the same order of $L^1(\lambda_2)$ error $(\log n)/n$. Remark
that for this frequency level $J_n$, we have to estimate $n/\log n$
coefficients. Estimations of one of these coefficients
$a_{\mathbf{k}}^{(J_n)}$, where $\mathbf{k}\in I_2(J_n)$, depends on
the number of observation in the square $\cI_{\mathbf{k}}^{(J_n)}$
associated this coefficient. The probability that no observation
"falls" in $\cI_{\mathbf{k}}^{(J_n)}$ is smaller than $n^{-1}$.
Thus, number of coefficients estimated with no observations is small
compare to the order of approach $(\log n)/n$ and is taken into
account in the variance term. Now, the problem is about finding a
$L^1-$ball of prediction rules such that for any integer $n$ the
approximation function $f_n$ belongs to such a ball. This problem
depends on the geometry of the border set $\partial A$. It arises
naturally since we chose a particular geometry for our partition:
dyadic partitions of the space $[0,1]^d$, and we have to pay a price
for this choice which has been made independently of the type of
functions to estimate. But this choice of geometry in our case is
the same as the one met in density approximation using approximation
theory while choosing a particular wavelet basis. Depending on the
type of Bayes rules we have to estimate, a special partition can be
considered. For example our "dyadic approach" is very well adapted
for the estimation of Bayes rules associated to chessboard (with the
value $1$ for black square and $-1$ for white square). This kind of
Bayes rules are very bad estimated by classification procedure
estimating the border since most of these procedure have regularity
assumptions which are not fulfilled in the case of chessboard.

We can extend our approach in several different ways. Consider the
dyadic partition of $[0,1]^d$ with frequency $J_n$. Instead of
choosing $1$ or $-1$ for each square of this partition (like in our
approach), we can do a least square regression in each cell of the
partition. Inside a square $Sq=\cI_{\mathbf{k}}^{(J_n)}$, where
$\mathbf{k}\in I_2(J_n)$, we can compute the line minimizing
$$\sum_{i=1}^n\1_{(2f(X_i)-1\neq Y_i,X_i\in Sq)},$$ where $f$ is
taken in the set of all indicators of half spaces of $[0,1]^d$
intersecting $Sq$. Of course, depending on the number of
observations inside the cell $Sq$ we can consider bigger classes of
functions than the one made of the indicators of half spaces. Our
classifier is close to the histogram estimator in density or
regression framework, which has been extend to smoother procedure.
The other way to extend our approach deals with the problem of the
underlying choice of geometry by taking $\cS$ for fundamental
system. One possible solution is to consider classifiers "adaptive
to the geometry". Using an adaptive procedure, for instance
aggregation procedure (cf. \cite{lec:05}), we can construct
classifiers adaptive to the "rotation" and "translation". Consider
the dyadic partition of $[0,1]^2$ at the frequency level $J_n$. We
can construct classifiers using the same procedure as (\ref{esti})
but for partitions obtained by translation of the dyadic partition
by $(n_1/(2^{J_n}\log n), n_2/(2^{J_n}\log n))$, where
$n_1,n_2=0,\ldots,\lceil \log n \rceil$. We can do the same thing by
aggregating classifiers obtained by the procedure (\ref{esti}) for
partitions obtained by rotation of center $(1/2,1/2)$ with angle
$n_3\pi/(2\log n)$, where $n_3=0,\ldots,\lceil \log n \rceil$, of
the initial dyadic partition. In this heuristic we don't discuss
about the way to solve problems near the border of $[0,1]^2$.

\par

\section{Proofs}

{\bf Proof of Theorem \ref{theoopenas}:} Since $\{\eta\geq1/2\}$ is
almost everywhere open there exists an open subset $\cO$ of
$[0,1]^d$ such that $\lambda_d(\{\eta\geq1/2\}\Delta\cO)=0$. If
$\cO$ is the empty set then take $g=-1$, otherwise, for all
$x\in\cO$ denote by $\cI_x$ the biggest subset
$\cI_{\mathbf{k}}^{(j)}$ for $j\in\mathbb{N}$ and $\mathbf{k}\in
I_d(j)$ such that $x\in\cI_{\mathbf{k}}^{(j)}$ and
$\cI_{\mathbf{k}}^{(j)}\subseteq\cO.$ Remark that $\cI_x$ exists
because $\cO$ is open. We can see that for any $y\in\cI_x$ we have
$\cI_y=\cI_x$, thus, $(\cI_x:x\in\cO)$ is a partition of $\cO$. We
denote by $I_\cO$ a subset of index $(j,\mathbf{k})$, where
$j\in\mathbb{N}, \mathbf{k}\in I_d(j)$ such that
$\{\cO_x:x\in\cO\}=\{\cI_\mathbf{k}^{(j)}:(j,\mathbf{k})\in
I_\cO\}.$ For any $(j,\mathbf{k})\in I_\cO$ we take
$a_\mathbf{k}^{(j)}=1$.

Take $\cO_1$ an open subset $\lambda_d$-almost everywhere equal to
$\{\eta<1/2\}$. If $\cO_1$ is the empty set then take $g=1$.
Otherwise, consider the set of index $I_{\cO_1}$ built in the same
way as previously, and for any $(j,\mathbf{k})\in I_{\cO_1}$ we take
$a_\mathbf{k}^{(j)}=-1$.

For all $(j,\mathbf{k})\notin I_\cO \cup I_{\cO_1}$, we take
$a_\mathbf{k}^{(j)}=0$. Consider
$$g=\sum_{j=0}^{+\infty}\sum_{\mathbf{k}\in I_d(j)}a_{\mathbf{k}}^{(j)}
\phi_{\mathbf{k}}^{(j)}.$$ It is easy to check that the function $g$
belongs to $\cF^{(d)}$ and satisfies the writing convention (W) and
that, for $\lambda_d-$almost $x\in[0,1]^d$, $g(x)=f_\eta(x)$.

{\bf{ Proof of Proposition \ref{p1}}:} Assume that
$\cF_w^{(d)}\neq \{\1_{[0,1]^d}\}$. Take $f\in\cF_w^{(d)}-
\{\1_{[0,1]^d}\}$. Consider the writing of $f$ in the system $\cS$
using the convention (W),
$$f=\sum_{j\in\mathbb{N}}\sum_{\mathbf{k}\in I_d(j)}a_{\mathbf{k}}^{(j)}
\phi_{\mathbf{k}}^{(j)},$$ where $a_{\mathbf{k}}^{(j)}\in\{-1,0,1\}$
for any $j\in\mathbb{N},\mathbf{k}\in I_d(j)$. Consider
$b_{\mathbf{k}}^{(j)}=|a_{\mathbf{k}}^{(j)}|$ for any
$j\in\mathbb{N},\mathbf{k}\in I_d(j)$. Take
$f_2=\sum_{j\in\mathbb{N}}\sum_{\mathbf{k}\in
I_d(j)}b_{\mathbf{k}}^{(j)}\phi_{\mathbf{k}}^{(j)}$. Remark that the
function $f_2\in\cF^{(d)}$ does not satisfy the writing convention
(W). We have $f_2=\1_{[0,1]^d}$. For any $j\in\mathbb{N}$ we have
\begin{equation}\label{equationcard}{\rm
card}\left\{\mathbf{k}\in I_d(j) :b_{\mathbf{k}}^{(j)}\neq0\right\}=
{\rm card}\left\{\mathbf{k}\in
I_d(j):a_{\mathbf{k}}^{(j)}\neq0\right\}.\end{equation}Moreover, one
coefficient $b_{\mathbf{k}}^{(j)}\neq0$ contributes to fill a cell
of Lebesgue measure $2^{-dj}$ among the hypercube $[0,1]^d$. Since
the mass total of $[0,1]^d$ is $1$, we have
\begin{equation}\label{equationmass}1=\sum_{j\in\mathbb{N}}\sum_{\mathbf{k}\in I_d(j)}
2^{-dj}{\rm card}\left\{\mathbf{k}\in
I_d(j):b_{\mathbf{k}}^{(j)}\neq0\right\}.\end{equation} Moreover,
$f\in\cF^{(d)}$ thus, for any $j\in\mathbb{N}$,
$$\lfloor w(j)\rfloor\geq{\rm
card}\left\{\mathbf{k}\in
I_d(j):a_{\mathbf{k}}^{(j)}\neq0\right\}.$$ We obtain the second
assertion of Proposition \ref{p1} by using the last inequality and
the both assertions (\ref{equationcard}) and (\ref{equationmass}).

Assume that $ \sum_{j=1}^{+\infty} 2^{-dj}\lfloor w(j)
\rfloor\geq1.$ For any integer $j\neq0$, denote by $\cI(j)$ the set
of indexes $\left\{(j,\mathbf{k}): \mathbf{k}\in I_d(j)\right\}$.

We use the natural order of $\mathbb{N}^{d+1}$ to order sets of
indexes. Take $\cI_w(1)$ the family of the first $\lfloor w(1)
\rfloor$ elements of $\cI(1)$. Denote by $\cI_w(2)$ the family made
of the first $\lfloor w(1) \rfloor$ elements of $\cI(1)$ and add, at
the end of this family in the correct order, the first $\lfloor w(2)
\rfloor$ elements $(2,\mathbf{k})$ of $\cI(2)$ such that
$\phi_{\mathbf{k}'}^{(1)}\phi_{\mathbf{k}}^{(2)}=0$ for any
$(1,\mathbf{k}')\in\cI_w(1)$,..., for the step $j$, construct the
family $\cI_w(j)$ made of all the elements of $\cI_w(j-1)$ in the
same order and add at the end of this family the indexes
$(j,\mathbf{k})$ in $\cI(j)$ among the first $\lfloor w(j) \rfloor$
elements of $\cI(j)$ such that
$\phi_{\mathbf{k}'}^{(J)}\phi_\mathbf{k}^{(j)}=0$ for any
$(J,\mathbf{k}')\in\cI_w(j-1)$. If there is no more index satisfying
this condition then we stop the construction otherwise we go on.
Denote by $\cI$ the final family obtained by this construction
($\cI$ may be finite or infinite). Then, we enumerate the indexes of
$\cI$ by $(j_1,\mathbf{k}_1)\prec(j_2,\mathbf{k}_2)\prec\cdots$. For
the first $(j_1,\mathbf{k}_1)\in\cI$ take
$a_{\mathbf{k}_1}^{(j_1)}=1$, for the second element
$(j_2,\mathbf{k}_2)\in\cI$ take $a_{\mathbf{k}_2}^{(j_2)}=-1$,etc. .
Consider the function
$$f=\sum_{j\in\mathbb{N}}\sum_{\mathbf{k}\in I_d(j)}a_{\mathbf{k}}^{(j)}
\phi_{\mathbf{k}}^{(j)}.$$ If the construction stops at a given
iteration $N$ then $f$ takes its values in $\{-1,1\}$ and the
writing convention (W) is fulfilled since every cells
$\cI_{\mathbf{k}}^{(j)}$ such that $a_{\mathbf{k}}^{(j)}\neq0$ has a
neighboring cell associated to a coefficient non equals to $0$ with
an opposite value. Otherwise, for any integer $j\neq0$, the number
of coefficient $a_{\mathbf{k}}^{(j)}$, for $\mathbf{k}\in I_d(j)$,
non equals to $0$ is $\lfloor w(j) \rfloor$ and the total mass of
cells $\cI_{\mathbf{k}}^{(j)}$ such that $a_{\mathbf{k}}^{(j)}\neq0$
is $\sum_{j\in\mathbb{N}}\sum_{\mathbf{k}\in I_d(j)}2^{-dj}{\rm
card}\left\{\mathbf{k}\in I_d(j):a_{\mathbf{k}}^{(j)}\neq0\right\}$
which is greater or equal to $1$ by assumption. Thus, all the
hypercube is filled by cells associated to coefficients non equal to
$0$. So $f$ takes its values in $\{-1,1\}$ and the writing
convention (W) is fulfilled since every cells
$\cI_{\mathbf{k}}^{(j)}$ such that $a_{\mathbf{k}}^{(j)}\neq0$ has a
neighboring cell associated to a coefficient non equals to $0$ with
an opposite value. Moreover $f\neq\1_{[0,1]^d}$.

{\bf{Proof of Theorem \ref{approxtheo}.}} Let $\pi=(P^X,\eta)$ be
a probability measure on $\cX\times\{-1,1\}$ belonging to
$\cP_{w,A}$. Denote by $f^*$ a Bayes classifier associated to
$\pi$ (for example $f^*={\rm{sign}}(2\eta-1)$) . We
have$$d_\pi(f,f^*)=(1/2)\mathbb{E}\left[|2\eta(X)-1||f(X)-f^*(X)|
\right]\leq (A/2)||f-f^*||_{L^1(\lambda_d)}.$$

Let $\epsilon>0$. Define by $J_\epsilon$ the smallest integer
satisfying
$$\sum_{j=J_\epsilon+1}^{+\infty}2^{-dj}\lfloor w(j) \rfloor<\frac{\epsilon}{A}.$$
We write $f^*$ in the fundamental system
$(\phi_{\mathbf{k}}^{(j)},j\in\mathbb{N},\mathbf{k}\in I_d(j))$
using the convention of writing of
section~\ref{subsectionapproximationresult} but we start at the
level of frequency $J_\epsilon$:
$$f^*=\sum_{\mathbf{k}\in I_d(J_\epsilon)}
A_{\mathbf{k}}^{(J_\epsilon)}\phi_{\mathbf{k}}^{(J_\epsilon)}
+\sum_{j=J_\epsilon+1}^{+\infty} \sum_{\mathbf{k}\in
I_d(j)}a_{\mathbf{k}}^{(j)}\phi_{\mathbf{k}}^{(j)}.$$ We consider
\begin{equation}\label{funcapprox}f_\epsilon=\sum_{\mathbf{k}\in I_d(J_\epsilon)}
B_{\mathbf{k}}^{(J_\epsilon)}\phi_{\mathbf{k}}^{(J_\epsilon)},\end{equation}
where
\begin{equation}\label{coefB}B_{\mathbf{k}}^{(J_\epsilon)}=\left\{\begin{array}{cl}
1 & \mbox{if }
p_{\mathbf{k}}^{(J_\epsilon)}>1/2\\
-1 & \mbox{otherwise }\end{array} \right.\end{equation} and
\begin{equation}\label{para} p_{\mathbf{k}}^{(J_\epsilon)}=\mathbb{P}(Y=1|X\in
I_{\mathbf{k}}^{(J_\epsilon)})=
\int_{I_{\mathbf{k}}^{(J_\epsilon)}}\eta(x)\frac{dP^X(x)}{P^X(I_{\mathbf{k}}^{(J_\epsilon)})},\end{equation}
for all $\mathbf{k}\in I_d(J_\epsilon)$. Note that, if
$A_{\mathbf{k}}^{(J_\epsilon)}\neq0$ then
$A_{\mathbf{k}}^{(J_\epsilon)}=B_{\mathbf{k}}^{(J_\epsilon)}$,
moreover $f^*$ take its values in $\{-1,1\}$, thus ,we have

\begin{eqnarray*}
||f_\epsilon-f^*||_{L^1(\lambda_d)} & = &
\sum_{\substack{\mathbf{k}\in
I_d(J_\epsilon)\\A_{\mathbf{k}}^{(J_\epsilon)}\neq0}}
\int_{I_{\mathbf{k}}^{(J_\epsilon)}} |f^*(x)-f_\epsilon(x)|dx+
\sum_{\substack{\mathbf{k}\in
I_d(J_\epsilon)\\A_{\mathbf{k}}^{(J_\epsilon)}=0}}
\int_{I_{\mathbf{k}}^{(J_\epsilon)}} |f^*(x)-f_\epsilon(x)|dx\\ &
\leq & 2^{-dJ_\epsilon+1}{\rm{card}}\left\{\mathbf{k}\in
I_d(J_\epsilon): A_{\mathbf{k}}^{(J_\epsilon)}=0\right\} \leq
2\sum_{j=J_\epsilon+1}^{+\infty}2^{-dj}\lfloor
w(j)\rfloor<2\epsilon/A.
\end{eqnarray*}

{\bf{Proof of Theorem \ref{estitheo}.}} Let $\pi=(P^X,\eta)$ be a
probability measure on $\cX\times\{-1,1\}$ satisfying (A1), (SMA)
and such that $f^*={\rm{sign}}(2\eta-1)$, a Bayes classifier
associated to $\pi$, belongs to  $\cF_w^{(d)}$ (a $L^1-$ball of
Bayes rules).

Let $\epsilon>0$ and $J_\epsilon$ the smallest integer satisfying
$\sum_{j=J_\epsilon+1}^{+\infty}2^{-dj}\lfloor
w(j)\rfloor<\epsilon/A$. We decompose the risk in the bias term and
variance term:
$$\cE(\hfn^{(J_\epsilon)})=\mathbb{E}\left[d_\pi(\hfn^{(J_\epsilon)},f^*) \right]\leq
\mathbb{E}\left[d_\pi(\hfn^{(J_\epsilon)},f_\epsilon)\right]+d_\pi(f_\epsilon,f^*),$$
where $\hfn^{(J_\epsilon)}$ is introduced in (\ref{esti}) and
$f_\epsilon$ in (\ref{funcapprox}).

Using the definition of $J_\epsilon$ and according to the
approximation Theorem (Theorem~\ref{approxtheory}), the bias term
satisfies:
$$d_\pi(f_\epsilon,f^*)\leq \epsilon.$$

For the variance term we have (using the notations introduced in
(\ref{esti}) and (\ref{coefB})):
\begin{eqnarray*}
\mathbb{E}\left[d_\pi(\hfn^{(J_\epsilon)},f_\epsilon)\right] & = &
\frac{1}{2}\left|\mathbb{E}\left[
Y(f_\epsilon(X)-\hfn^{(J_\epsilon)}(X))\right]\right|\leq
\frac{1}{2}\mathbb{E}\left[ \int_{[0,1]^d}|f_\epsilon(x)
-\hfn^{(J_\epsilon)}(x)|dP^X(x)\right]\\
& = & \frac{1}{2}\sum_{\mathbf{k}\in
I_d(J_\epsilon)}\mathbb{E}\left[
\int_{I_{\mathbf{k}}^{(J_\epsilon)}}|B_{\mathbf{k}}^{(J_\epsilon)}
-\hat{A}_{\mathbf{k}}^{(J_\epsilon)}|
dP^X\right]\\
& \leq & \frac{A}{2^{dJ_\epsilon+1}}\sum_{\mathbf{k}\in
I_d(J_\epsilon)}
\mathbb{E}[|B_{\mathbf{k}}^{(J_\epsilon)}-\hat{A}_{\mathbf{k}}^{(J_\epsilon)}|]
  \leq  \frac{A}{2^{dJ_\epsilon}}\sum_{\mathbf{k}\in I_d(J_\epsilon)}
 \mathbb{P}\left(|B_{\mathbf{k}}^{(J_\epsilon)}-\hat{A}_{\mathbf{k}}^{(J_\epsilon)}|=2
\right).
\end{eqnarray*}
Let $\mathbf{k}\in I_d(J_\epsilon)$. For any $m\in\{0,\ldots,n\}$,
we introduce the sets
$$\Omega_{\mathbf{k}}^{(m)}=\left\{{\rm Card}\{ i\in\{1,\ldots,n\}: X_i\in
I_{\mathbf{k}}^{(J_\epsilon)}\}=m\right\}$$ and
$$ \Omega_{\mathbf{k}}=\left\{ {\rm{card}}
\left\{i\in\{1,\ldots,n\}:\begin{array}{l}X_i\in
I_{\mathbf{k}}^{(J_\epsilon)},\\ Y_i=1\end{array} \right\}\leq
{\rm{card}} \left\{i\in\{1,\ldots,n\}:\begin{array}{l}X_i\in
I_{\mathbf{k}}^{(J_\epsilon)},\\ Y_i=-1\end{array} \right\}
\right\}.$$ We have
\begin{equation*}
\mathbb{P}(\hat{A}_{\mathbf{k}}^{(J_\epsilon)}=-1)  =
\mathbb{P}(\Omega_{\mathbf{k}}^{(0)c}\cap\Omega_{\mathbf{k}})
+\mathbb{P}(\Omega_{\mathbf{k}}^{(0)})
\end{equation*}and
$$
  \mathbb{P}(\Omega_{\mathbf{k}}^{(0)c}\cap\Omega_{\mathbf{k}})
=\sum_{m=1}^n
  \mathbb{P}(\Omega_{\mathbf{k}}^{(m)}\cap\Omega_{\mathbf{k}})
  = \sum_{m=1}^n
  \mathbb{P}(\Omega_{\mathbf{k}}|\Omega_{\mathbf{k}}^{(m)})
  \mathbb{P}(\Omega_{\mathbf{k}}^{(m)}).
$$

Moreover, denote by $Z_1,\ldots,Z_n$ some variables i.i.d. with a
Bernoulli with parameter $p_{\mathbf{k}}^{(J_\epsilon)}$ for common
probability distribution ($p_{\mathbf{k}}^{(J_\epsilon)}$ is
introduced in (\ref{para}) and is equal to $\mathbb{P}(Y=1|X\in
I_{\mathbf{k}}^{(J_\epsilon)})$), we have for any $m=1,\ldots,n$,
$$\mathbb{P}(\Omega_{\mathbf{k}}|\Omega_{\mathbf{k}}^{(m)})=
\mathbb{P}\left(\frac{1}{m}\sum_{i=1}^m Z_i\leq \frac{1}{2}
\right).$$ Concentration inequality of Hoeffding leads to
\begin{equation}\label{hoeffding}
\mathbb{P}\left(\frac{1}{m}\sum_{i=1}^m Z_i\geq
p_{\mathbf{k}}^{(J_\epsilon)}+t \right)\leq \exp(-2mt^2) \mbox{ and
} \mathbb{P}\left(\frac{1}{m}\sum_{i=1}^m Z_i\leq
p_{\mathbf{k}}^{(J_\epsilon)}-t \right)\leq \exp(-2mt^2),
\end{equation}
for all $t>0$ and $m=1,\ldots,n$.

Denote by $a_{\mathbf{k}}^{(J_\epsilon)}$ the probability
$\mathbb{P}\left(X\in I_{\mathbf{k}}^{(J_\epsilon)} \right)$. If
$p_{\mathbf{k}}^{(J_\epsilon)}>1/2$, applying second inequality of
(\ref{hoeffding}) leads to
\begin{eqnarray*}
\lefteqn{\mathbb{P}\left(|B_{\mathbf{k}}^{(J_\epsilon)}-\hat{A}_{\mathbf{k}}^{(J_\epsilon)}|=2
\right)  =  \mathbb{P}(\hat{A}_{\mathbf{k}}^{(J_\epsilon)}=-1)}\\
& \leq & \sum_{m=1}^n\mathbb{P}\left[\frac{1}{m}\sum_{j=1}^m Z_j\leq
p_{\mathbf{k}}^{(J_\epsilon)}-(p_{\mathbf{k}}^{(J_\epsilon)}- 1/2)
\right]\left(\begin{array}{c} n\\m
\end{array}\right)(a_{\mathbf{k}}^{(J_\epsilon)})^m(1-a_{\mathbf{k}}^{(J_\epsilon)})^{n-m}\\
&+&\mathbb{P}(\Omega_{\mathbf{k}}^{(0)})\\
& \leq &
\sum_{m=0}^n\exp\left(-2m(p_{\mathbf{k}}^{(J_\epsilon)}-1/2)^2
\right)\left(\begin{array}{c} n\\m
\end{array}\right)(a_{\mathbf{k}}^{(J_\epsilon)})^m(1-a_{\mathbf{k}}^{(J_\epsilon)})^{n-m}\\
&=&\left(1-
a_{\mathbf{k}}^{(J_\epsilon)}(1-\exp(-2(p_{\mathbf{k}}^{(J_\epsilon)}-1/2)^2))\right)^n\\
&\leq&
\exp\left(-na(1-\exp(-2(p_{\mathbf{k}}^{(J_\epsilon)}-1/2)^2))2^{-dJ_\epsilon}
\right).
\end{eqnarray*}

If $p_{\mathbf{k}}^{(J_\epsilon)}<1/2$ then similar arguments used
in the previous case and first inequality of (\ref{hoeffding}) lead
to
\begin{eqnarray*}\mathbb{P}\left(|B_{\mathbf{k}}^{(J_\epsilon)}-\hat{A}_{\mathbf{k}}^{(J_\epsilon)}|=2
\right)& =&
\mathbb{P}(\hat{A}_{\mathbf{k}}^{(J_\epsilon)}=1)\\&\leq&
\exp\left(-na(1-\exp(-2(p_{\mathbf{k}}^{(J_\epsilon)}-1/2)^2))2^{-dJ_\epsilon}
\right).\end{eqnarray*} If $p_{\mathbf{k}}^{(J_\epsilon)}=1/2$, we
use
$\mathbb{P}\left(|B_{\mathbf{k}}^{(J_\epsilon)}-\hat{A}_{\mathbf{k}}^{(J_\epsilon)}|=2
\right)\leq 1$. Like in the proof of Theorem \ref{approxtheo}, we
use the writing $$f^*=\sum_{\mathbf{k}\in I_d(J_\epsilon)}
A_{\mathbf{k}}^{(J_\epsilon)}\phi_{\mathbf{k}}^{(J_\epsilon)}
+\sum_{j=J_\epsilon+1}^{+\infty} \sum_{\mathbf{k}\in
I_d(j)}a_{\mathbf{k}}^{(j)}\phi_{\mathbf{k}}^{(j)}.$$ Since
$P^X(\eta=1/2)=0$, if $A_{\mathbf{k}}^{(J_\epsilon)}\neq 0$ then
$p_{\mathbf{k}}^{(J_\epsilon)}\neq 1/2$. Thus, the variance term
satisfies:
\begin{eqnarray*}
\lefteqn{\mathbb{E}\left[d_\pi(\hfn,f_\epsilon^*)\right]}\\ & \leq &
\frac{A}{2^{dJ_\epsilon}}\left( \sum_{\substack{\mathbf{k}\in
I_d(J_\epsilon)\\A_{\mathbf{k}}^{(J_\epsilon)}\neq 0}}
\mathbb{P}\left(|B_{\mathbf{k}}^{(J_\epsilon)}-\hat{A}_{\mathbf{k}}^{(J_\epsilon)}|=2
\right)+ \sum_{\substack{\mathbf{k}\in
I_d(J_\epsilon)\\A_{\mathbf{k}}^{(J_\epsilon)}= 0}}
\mathbb{P}\left(|B_{\mathbf{k}}^{(J_\epsilon)}-\hat{A}_{\mathbf{k}}^{(J_\epsilon)}|=2
\right) \right)\\
& \leq & \frac{A}{2^{dJ_\epsilon}} \sum_{\substack{\mathbf{k}\in
I_d(J_\epsilon)\\A_{\mathbf{k}}^{(J_\epsilon)}\neq
0}}\exp\left(-na(1-\exp(-2(p_{\mathbf{k}}^{(J_\epsilon)}-1/2)^2))2^{-dJ_\epsilon}
\right)+A\epsilon.
\end{eqnarray*}
If $A_{\mathbf{k}}^{(J_\epsilon)}\neq 0$ then $\eta>1/2$ or
$\eta<1/2$ over the whole set $I_{\mathbf{k}}^{(J_\epsilon)}$, so
$$\left\vert
\frac{1}{2}-p_{\mathbf{k}}^{(J_\epsilon)}\right\vert=\int_{I_{\mathbf{k}}^{(J_\epsilon)}}
\left|
\eta(x)-\frac{1}{2}\right|\frac{dP^X(x)}{P^X(I_{\mathbf{k}}^{(J_\epsilon)})}.$$
Moreover $\pi$ satisfies $\mathbb{P}\left(|2\eta(X)-1|\geq
h\right)=1$, so
$$\left\vert
\frac{1}{2}-p_{\mathbf{k}}^{(J_\epsilon)}\right\vert\geq
\frac{h}{2}.$$

We have shown that for all $\epsilon>0$,
$$\cE(\hfn)=\mathbb{E}\left[d_\pi(\hfn,f^*) \right]\leq
(1+A)\epsilon+\exp\left(-na(1-\exp(-2(h/2)^2))2^{-dJ_\epsilon}
\right),$$ where $J_\epsilon$ is the smallest integer satisfying
$\sum_{j=J_\epsilon+1}^{+\infty}2^{-dj}\lfloor w(j)
\rfloor<\epsilon/A$.

{\bf{Proof of Theorem \ref{optimal}}.} For all $q\in\mathbb{N}$ we
consider $G_q$ a net of $[0,1]^d$ defined
by:$$G_q=\left\{\left(\frac{2k_1+1}{2^{q+1}},\ldots,\frac{2k_d+1}{2^{q+1}}
\right): (k_1,\ldots,k_d)\in\{0,\ldots,2^q-1\}^d \right\}$$ and the
function $\eta_q$ from $[0,1]^d$ to $G_q$ such that $\eta_q(x)$ is
the closest point of $G_q$ from $x$ (in the case of ex aequo, we
choose the smallest point for the usual order on $\mathbb{R}^d$).
Associated to this grid, the partition
${\cX'}_1^{(q)},\ldots,{\cX'}_{2^{dq}}^{(q)}$ of $[0,1]^d$ is
defined by $x,y\in{\cX'}_i^{(q)}$ iff $\eta_q(x)=\eta_q(y)$ and we
use a special indexation for this partition: denote by
${x'}_{k_1,\ldots,k_d}^{(q)}=\left(\frac{2k_1+1}{2^{q+1}},\ldots,\frac{2k_d+1}{2^{q+1}}
\right)$ and we say that ${x'}_{k_1,\ldots,k_d}^{(q)}\prec
{x'}_{{k'}_1,\ldots,{k'}_d}^{(q)}$ if $$
\eta_{q-1}({x'}_{k_1,\ldots,k_d}^{(q)})\prec
\eta_{q-1}({x'}_{{k'}_1,\ldots,{k'}_d}^{(q)})$$ or
$$\eta_{q-1}({x'}_{k_1,\ldots,k_d}^{(q)})=
\eta_{q-1}({x'}_{{k'}_1,\ldots,{k'}_d}^{(q)})\mbox{ and }
(k_1,\ldots,k_d)<({k'}_1,\ldots,{k'}_d),$$ for the usual order on
$\mathbb{N}^d$. Thus, the partition
$({\cX'}_j^{(q)}:j=1,\ldots,2^{dq})$ has an increasing indexation
according to the order of $({x'}^{(q)}_{k_1,\ldots,k_d})$ for the
order defined above. This order take care of the previous
partition by splitting blocks in the right given order and inside
a block of a partition we take the natural order of
$\mathbb{N}^d$. We introduce an other parameter
$m\in\{1,\ldots,2^{qd}\}$ and we define for all $i=1,\ldots,m$,
$\cX_i^{(q)}={\cX'}_i^{(q)}$ and $\cX_0^{(q)}=[0,1]^d-\cup_{i=1}^m
\cX_i^{(q)}$. Parameters $q$ and $m$ will be chosen later. We
consider $W\in[0,m^{-1}]$, chosen later, and define the function
$f_X$ from $[0,1]^d$ to $\mathbb{R}$ by $f_X=W/\lambda_d(\cX_1)$
(where $\lambda_d$ is the Lebesgue measure on $[0,1]^d$) on
$\cX_1,\ldots,\cX_m$ and $(1-mW)/\lambda_d(\cX_0)$ on $\cX_0$. We
denote by $P^X$ the probability distribution on $[0,1]^d$ with the
density $f_X$ w.r.t. the Lebesgue measure. For all
$\sigma=(\sigma_1,\ldots,\sigma_m)\in\Omega=\{-1,1\}^m$ we
consider $\eta_\sigma$ defined for any  $x\in[0,1]^d$ by
\begin{equation*}
\eta_\sigma(x)=\left\{
  \begin{array}{cc}
  \frac{1+\sigma_j h}{2} & \mbox{ if } x\in\cX_j, j=1,\ldots,m,\\
  1 & \mbox{ if } x\in\cX_0.
  \end{array}\right.
\end{equation*}We have a set of probability measures $\{\pi_\sigma:\sigma\in\Omega\}$
on $[0,1]^d\times\{-1,1\}$ indexed by the hypercube $\Omega$ where
$P^X$ is the marginal on $[0,1]^d$ of $\pi_\sigma$ and
$\eta_\sigma$ its conditional probability function of $Y=1$ given
$X$. We denote by $f_\sigma^*$ the Bayes rule associated to
$\pi_\sigma$, we have $f_\sigma^*(x)=\sigma_j$ if $x\in\cX_j$ for
$j=1,\ldots,m$ and $1$ if $x\in\cX_0$, for any $\sigma\in\Omega$.

Now we give conditions on $q,m$ and $W$ such that for all $\sigma$
in $\Omega$, $\pi_\sigma$ belongs to ${\cP}_{w,h,a,A}$. If we take
\begin{equation}\label{condw}W=2^{-dq},\end{equation} then $P^X<<\lambda$ and $\forall x\in[0,1]^d,
a\leq dP^X/d\lambda(x) \leq A$. We have clearly $\vert
2\eta(x)-1\vert\geq h$ for any  $x\in[0,1]^d$. We can see that
$f_\sigma^*\in{\cF}_w^{(d)}$ for all $\sigma\in\{-1,1\}^m$ iff
$$\begin{array}{l}
\lfloor w (q+1)\rfloor\geq \inf(x\in2^d\mathbb{N}:x\geq m)\\
\lfloor w(q) \rfloor\geq\left\{\begin{array}{lll}2^d-1 &
\mbox{if } & m<2^d\\
\inf(x\in2^d\mathbb{N}:x\geq 2^{-d}m) & \mbox{otherwise} & \end{array}\right.\\
\ldots\\
\lfloor w(1)\rfloor\geq\left\{\begin{array}{lll}2^d-1 &
\mbox{if } & m<2^{dq}\\
\inf(x\in2^d\mathbb{N}:x\geq 2^{-dq}m) & \mbox{otherwise} & \end{array}\right.\\
\lfloor w(0) \rfloor\geq 1
\end{array}.$$
Since we have $\lfloor w  (j)\rfloor\geq 2^d-1$ for all $j\geq1$
and $\lfloor w(0) \rfloor=1$, and $\lfloor w(j-1)\rfloor\geq
\lfloor w(j)\rfloor/2^d$, then $f_\sigma^*\in{\cF}_w^{(d)}$ for
all $\sigma\in\Omega$ iff
\begin{equation}\label{condmq} \lfloor w(q+1) \rfloor\geq
\inf(x\in2^d\mathbb{N}:x\geq m).\end{equation}

Take $q,m$ and $W$ such that (\ref{condw}) and (\ref{condmq}) are
fulfilled then, $\{\pi_\sigma:\sigma\in \Omega\}$ is a subset of
${\cP}_{w,h,a,A}$. Let $\sigma\in\Omega$ and $\hat{f}_n$ be a
classifier, we have
\begin{eqnarray*}
\mathbb{E}_{\pi_\sigma}\left[ R(\hat{f}_n)-R^*\right] & = &
(1/2)\mathbb{E}_{\pi_\sigma}\left[ \vert2\eta_\sigma(X)-1 \vert
\vert \hat{f}_n(X)-f^*_\sigma(X)\vert \right]\\ & \geq &
(h/2)\mathbb{E}_{\pi_\sigma}\left[ \vert
\hat{f}_n(X)-f^*_\sigma(X)\vert \right]\\ & \geq &
(h/2)\mathbb{E}_{\pi_\sigma}\left[\sum_{i=1}^m \int_{\cX_i}\vert
\hat{f}_n(x)-f^*_\sigma(x) \vert dP^X(x)+\int_{\cX_0}\vert
\hat{f}_n(x)-f^*_\sigma(x) \vert dP^X(x) \right]\\ & \geq &
(Wh/2)\sum_{i=1}^m \mathbb{E}_{\pi_\sigma}\left[ \int_{\cX_i}\vert
\hat{f}_n(x)-\sigma_i \vert\frac{dx}{\lambda(\cX_1)}\right]\\
& \geq & (Wh/2)\mathbb{E}_{\pi_\sigma}\left[\sum_{i=1}^m\left\vert
\sigma_i-\int_{\cX_i}\hat{f}_n(x)\frac{dx}{\lambda(\cX_1)}\right\vert\right].
\end{eqnarray*}
We deduce that $$\inf_{\hat{f}_n} \sup_{\pi\in{\cP}_{w,h,a,A}}
 \cE_{\pi}(\hat{f}_n)\geq (Wh/2)
\inf_{\hat{\sigma}_n\in[-1,1]^m} \sup_{\sigma\in\{-1,1\}^m}
\mathbb{E}_{\pi_\sigma}\left[ \sum_{i=1}^m \vert
\sigma_i-\hat{\sigma}_i\vert\right].$$

Now, we control the Hellinger distance between two neighbouring
probability measures. Let $\rho$ be the Hamming distance on
$\Omega$. Let $\sigma, \sigma'$ in $\Omega$ such that
$\rho(\sigma,\sigma')=1$. We have $$H^2(\pi_{\sigma}^{\otimes n},
\pi_{\sigma'}^{\otimes n})=2\left(1-\left(
1-\frac{H^2(\pi_\sigma,\pi_{\sigma'})}{2}\right)^n \right),$$ and
a straightforward calculus leads to
$H^2(\pi_\sigma,\pi_{\sigma'})=2W\left(1-\sqrt{1-h^2} \right)$.
Take \begin{equation}\label{coucou}W=1/n,\end{equation} thus, for
any integer $n$, we have $H^2(\pi_{\sigma}^{\otimes n},
\pi_{\sigma'}^{\otimes n})\leq \beta<2$ where
$\beta=2\left(1-\exp(1-\sqrt{1-h^2})\right)$. The Assouad's Lemma
(cf. \cite{lec3:06}) yields $\inf_{\hat{\sigma}_n\in[-1,1]^m}
\sup_{\sigma\in\{-1,1\}^m} \mathbb{E}_{\pi_\sigma}\left[
\sum_{i=1}^m \vert \sigma_i-\hat{\sigma}_i\vert\right]\geq
\frac{m}{4}\left(1-\frac{\beta}{2} \right)^2$. We conclude that
$$\inf_{\hat{f}_n} \sup_{\pi\in\cP_{w,h,a,A}}
 \cE_{\pi}(\hat{f}_n)\geq
Wh\frac{m}{8}\left(1-\frac{\beta}{2} \right)^2.$$

According to (\ref{condw}), (\ref{condmq}) and (\ref{coucou}) we
take $W=2^{-dq}=1/n, q=\left\lfloor \log n/(d\log 2)\right\rfloor,
m=\lfloor w\left(\lfloor \log
n/(d\log2)\rfloor+1\right)\rfloor-(2^d-1)$. For these values we
have
$$ \inf_{\hat{f}_n}\sup_{\pi\in\cP_{w,h,a,A}}
\cE_\pi(\hat{f}_n)\geq C_0 n^{-1}\left(\lfloor w\left(\lfloor \log
n/(d\log2)\rfloor+1\right)\rfloor-(2^d-1)\right).$$ where
$C_0=(h/8)\exp\left(-(1-\sqrt{1-h^2}) \right).$

{\bf Proof of Corollary \ref{theodgl}:} It suffices to apply
Theorem \ref{optimal} to the function $w$ defined by $w(j)=2^{dj}$
for any integer $j$ and $a=A=1$ for $P^X=\lambda_d$.

{\bf Proof of Theorem \ref{t2}:}
\begin{enumerate}

\item If we assume that $J_\epsilon\geq K$ then
$\sum_{j=J_\epsilon+1}^{+\infty}2^{-dj}\lfloor
w_{K}^{(d)}(j)\rfloor=(2^{dK})/(2^{dJ_\epsilon}(2^d-1))$. We take
$$J_\epsilon=\left\lceil \frac{\log
\left((A2^{dK})/(\epsilon(2^d-1)) \right)}{d\log 2}\right\rceil$$
and $\epsilon_n$ the unique solution of
$(1+A)\epsilon_n=\exp(-nC\epsilon_n)$, where
$C=a(1-e^{-h^2/2})(2^d-1)[A2^{d(K+1)}]^{-1}$. Thus,
$\epsilon_n\leq (\log n)/(Cn)$. For $J_n=J_{\epsilon_n}$, we have
$$\cE\left(\hat{f}_n^{(J_n)} \right)\leq  C_{K,d,h,a,A}\frac{\log n}{n},$$
for any integer $n$ such that $\log n\geq 2^{d(K+1)}(2^d-1)^{-1}$
and $J_n\geq K$, where $C_{K,d,h,a,A}=2(1+A)/C$.

If we have $\lfloor \log n/(d\log2)\rfloor\geq 2$ then $\lfloor
w\left(\lfloor \log
n/(d\log2)\rfloor+1\right)\rfloor-(2^d-1)\geq2^{d}$, so we obtain
the lower bound with the constant $C_{0,K}=2^d C_0$ and if
$\lfloor \log n/(d\log2)\rfloor\geq K$ the constant can be
$C_{0,K}=C_0(2^{dK}-(2^d-1))$.

\item If we have $J_\epsilon\geq N^{(d)}(\alpha)$, then $
\sum_{j=J_\epsilon+1}^{+\infty}2^{-dj}\lfloor
w_{\alpha}^{(d)}(j)\rfloor\leq
(2^{d(1-\alpha)J_\epsilon}(2^{d(1-\alpha)}-1))^{-1}$. We take
$$J_\epsilon=\left\lceil
\frac{\log(A/(\epsilon(2^{d(1-\alpha)}-1)))}{d(1-\alpha)\log
2}\right\rceil.$$ Denote by $\epsilon_n$ the unique solution of
$(1+A)\epsilon_n=\exp(-nC\epsilon_n^{1/(1-\alpha)})$ where
$C=a(1-e^{-h^2/2})2^{-d}(A^{-1}(2^{d(1-\alpha)}-1))^{1/(1-\alpha)}.$
We have $\epsilon_n\leq (\log n/(nC))^{1-\alpha}$. For
$J_n=J_{\epsilon_n}$, we have
\begin{equation*} \cE\left(\hat{f}_n^{(J_n)} \right)  \leq
\frac{2(1+A)A}{2^{d(1-\alpha)}-1}\left[\frac{2^d}{a(1-e^{-h^2/2})}
\right]^{1-\alpha}\left( \frac{\log
n}{n}\right)^{1-\alpha}.\end{equation*}

For the lower bound we have for any integer $n$,
$$ \inf_{\hat{f}_n}\sup_{\pi\in\cP_{\alpha}^{(d)}}
\cE_\pi(\hat{f}_n)\geq C_0\max \left(1,
n^{-1}\left(2^dn^\alpha-(2^d-1)\right)\right).$$
\end{enumerate}

{\bf Proof of Theorem \ref{circle}:} Let $\epsilon>0$. Denote by
$\epsilon_0$ the greatest positive number satisfying
$\delta(\epsilon_0)\epsilon_0^2\leq \epsilon$. Consider
$N(\epsilon_0)=\cN(\partial A,\epsilon_0,||.||_\infty)$ and
$x_1,\ldots,x_{N(\epsilon_0)}\in\mathbb{R}^2$ such that $\partial
A\subset \cup_{j=1}^{N(\epsilon_0)}B_{\infty}(x_j,\epsilon_0)$.
Since $2^{-J_{\epsilon_0}}\geq \epsilon_0$, only nine dyadic sets
of frequency $J_{\epsilon_0}$ can be used to cover a ball of
radius $\epsilon_0$ for the infinity norm of $\mathbb{R}^2$. Thus,
we only need $9N(\epsilon_0)$ dyadic sets of frequency
$J_{\epsilon_0}$ to cover $\partial A$. Consider the partition of
$[0,1]^2$ by dyadic sets of frequency $J_{\epsilon_0}$. Except on
the $9N(\epsilon_0)$ dyadic sets used to cover the border
$\partial A$, the prediction rule $f_A$ is constant, equal to $1$
or $-1$, on the other dyadic sets. Thus, by taking
$f_{\epsilon_0}=\sum_{k_1,k_2=0}^{2^{J_{\epsilon_0}}-1}a_{k_1,k_2}^{(J_{\epsilon_0})}\phi^{(J_{\epsilon_0})}_{k_1,k_2},$
where $a_{k_1,k_2}^{(J_{\epsilon_0})}$ is equal to one value of
$f_A$ in the dyadic set $\cI_{k_1,k_2}^{(J_{\epsilon_0})}$, we
have
$$||f_{\epsilon_0}-f_A||_{L^1(\lambda_2)}
\leq 9N(\epsilon_0)2^{-2J_{\epsilon_0}}\leq36\delta(\epsilon_0)\epsilon_0^2\leq36\epsilon.$$


\bibliographystyle{plainnat}
\footnotesize{}

\end{document}